\documentclass[letterpaper, 10 pt, conference,draftcls,onecolumn]{ieeeconf} 
\IEEEoverridecommandlockouts 
\usepackage[belowskip=-10pt,aboveskip=0pt]{caption}
\usepackage{graphicx}
\usepackage{subfigure}
\usepackage{amsmath}
\usepackage{amsfonts}
\usepackage[psamsfonts]{amssymb}
\usepackage{comment}
\usepackage{times}
\usepackage{changepage}
\usepackage{cases}
\usepackage{mathtools}
\usepackage{color}
\usepackage{booktabs}
\usepackage{array}

\newtheorem{definition}{Definition}

\newtheorem{lemma}{Lemma}

\newtheorem{theorem}{Theorem}

\IEEEtriggeratref{999}

\newcommand{\R}{\mathbb{R}}
\newcommand{\N}{\mathbb{N}}

\newcommand{\kde}{\kappa(\delta, \epsilon)}

\newcommand{\dom}{\mathcal{X}}
\newcommand{\Adj}{\textnormal{Adj}}
\renewcommand{\P}{\mathbb{P}}
\newcommand{\Q}{\mathcal{Q}}

\newcolumntype{M}{>{\centering\arraybackslash}m{\dimexpr.25\linewidth-2\tabcolsep}}

\title{\LARGE \bf
Differentially Private Cloud-Based Multi-Agent \\ Optimization with Constraints
}

\author{M.T. Hale and M. Egerstedt$^\dag$\thanks{$^\dag$The authors are with the School of Electrical and Computer Engineering, Georgia Institute of
Technology, Atlanta, GA 30332, USA. Email: \texttt{\{matthale, magnus\}@gatech.edu}. Research supported in
part by  the NSF under Grant CNS-1239225.}
}

\begin{document}

\maketitle
\thispagestyle{empty}
\pagestyle{empty}

\begin{abstract}
We present an optimization framework that solves
constrained multi-agent optimization problems
while keeping each agent's state differentially private.
 The agents
in the network seek to optimize a local objective
function in the presence of global constraints.
Agents communicate only through a trusted cloud computer
and the cloud also performs computations based on global
information.
The cloud computer modifies the results of such computations
before they are sent to the agents in order to
guarantee that the agents' states are kept private. 
We show that under
mild conditions each agent's optimization problem
converges in mean-square to
its unique solution while each agent's state
is kept differentially private.
A numerical simulation is provided
to demonstrate the viability of this approach. 
\end{abstract}

\section{Introduction}
Multi-agent optimization problems arise naturally in a
variety of settings including wireless sensor networks \cite{khan09,trigoni12},
robotics \cite{soltero13}, communications \cite{kelly98,chiang07}, and power systems \cite{caron10}.
A number of approaches to solving such problems have been proposed.
In \cite{nedich09}, a distributed approach was introduced
in which each agent relies only on local information in a time-varying network to
solve convex optimization problems with non-differentiable objective functions.
In \cite{wei13} the authors present a distributed implementation of 
Newton's method to solve network utility maximization problems.
In \cite{nedich10}, a distributed approach
to solving consensus and set-constrained optimization
problems over time-varying networks is proposed.

In the current paper,
we consider inequality constrained multi-agent
optimization problems and add the additional requirement that
each agent's state value be kept private.
In \cite{hale14} we proposed a cloud-based architecture
for multi-agent optimization 
(without privacy) 
in which a cloud computer is used
to handle communications and necessary global computations. 
We use this architecture
as our starting point, though we change the cloud's 
role substantially so that it keeps each agent's state
differentially private. 

The notion of differential privacy was
first established in the database literature as a way of providing
strong
practical guarantees of privacy to users contributing
personal data to a database \cite{dwork06a,dwork06b}. 
Essentially
differential privacy
guarantees that any query of a
database does not change by much if a single element
of that database changes or is deleted, thereby
concealing individual database entries. One appealing aspect
of differential privacy is that post-processing
cannot weaken its privacy guarantees, 
allowing for the results of differentially
private queries to be freely processed \cite{dwork13}.
In addition, differential privacy guarantees that an entry in
a database cannot be determined exactly even if a malicious
adversary has any arbitrary side information, e.g.,
other database entries \cite{kasiv08}.

Recently,
differential privacy was extended
to dynamical systems in \cite{pappas14}. 
Roughly, 
a system is differentially
private if input signals which are close in the 
input space produce output signals which are close
in the output space. This definition provides the same
resilience to post-processing and protection against
adversaries with arbitrary side-information mentioned
above. 
It is this
form of differential privacy which we apply
here to constrained multi-agent optimization. 

In the context of optimization,
differential privacy has been applied in a
number of ways. 
It was used to carry out optimization
of piecewise-affine functions in \cite{han14} to keep
certain terms in the objective functions private.  
In \cite{huang14}, the authors solve
a distributed optimization problem in which the
agents' objective functions must be kept private.
In \cite{hsu14}, differentially private linear programs
are solved while constraints or the objective
function are kept private. 
In the current paper
a saddle point finding algorithm
in the vein of \cite{hale14} is used to keep states private. When computations
are performed using this algorithm, noises are added
in accordance with the framework for
differentially private dynamic systems set forth in
\cite{pappas14} in order provide the guarantee of privacy.


The rest of the paper is organized as follows.
First, Section II reviews the relevant results
in the existing research literature.
Then Section III formulates the specific
problem to be solved here and proves that 
it can be solved privately.
Next, Section IV
provides simulation results to support the theoretical
developments made. Finally, Section V
concludes the paper.

\section{Review of Relevant Results} \label{sec:review}
\subsection{Problem Overview} \label{ss:overview}
Let there be a network of $n$ agents indexed by 
$i \in A$, $A = \{1, \ldots, n\}$. Let each agent have a
scalar state $x_i \in \R$ and let each agent have a
local, convex objective function $f_i : \R \to \R$
that is $C^2$ in $x_i$. 
The function $f_i$ is assumed to be private in the sense that agent $i$
does not share it with other agents or the cloud. 
Let the agents be subject to $m$ global inequality constraints, $g_j : \R^n \to \R$,
where we require
\begin{equation} \label{eq:gj}
g_j(x) \leq 0
\end{equation}
for all $j \in \{1, \ldots, m\}$.
Here $x = (x_1, \ldots, x_n)^T \in \R^n$ contains all states
in the network and we assume $g_j \in C^2$ is a convex function in $x$ for every $j$. 
We let $x$ takes values in a nonempty, compact, convex
set $\dom \subset \R^n$ and assume that Slater's condition holds, namely
that
there is some $\bar{x} \in \dom$ satisfying
$g(\bar{x}) < 0$. 

We can form an equivalent global optimization problem through 
defining a global objective function $f$ by
\begin{equation}
f(x) = \sum_{i=1}^{n} f_i(x_i).
\end{equation}
Forming the Lagrangian of the global optimization problem, we have
\begin{equation} \label{eq:bigl}
L(x, \mu) = f(x) + \mu^Tg(x),
\end{equation}
where $\mu$ is a vector of Kuhn-Tucker (KT)
multipliers in the non-negative orthant of
$\R^m$, denoted $\R^m_{+}$. By definition
$L(\cdot, \mu)$ is convex and $L(x, \cdot)$ 
is concave. 

Seminal work by Kuhn and Tucker \cite{uzawa58b} showed that constrained
optima of $f$ subject to $g$ are saddle points of $L$.
We assume that $f$ has a unique constrained minimum and that the constrained
minimum of $f$ is a regular point of $g$ so that there is a unique
saddle point of $L$ \cite{chachuat07}. The saddle point of $L$,
denoted $(\hat{x}, \hat{\mu})$, satisfies the inequalities
\begin{equation} \label{eq:lag}
L(\hat{x}, \mu) \leq L(\hat{x}, \hat{\mu}) \leq L(x, \hat{\mu})
\end{equation}
for all admissible $x$ and $\mu$. 
One method for finding saddle points from
some initial point $(x(0), \mu(0))$ is due to
Bakushinskii and Polyak \cite{bakushinskii74}. 
Letting $P_{\dom}[\cdot]$ denote the projection
onto the set $\dom$ and $[\cdot]_+$ the projection onto
$\R^m_+$, we can iteratively compute new values of $x$ and $\mu$
according to 
\begin{equation} \label{eq:bigx0}
x(k+1) = P_{\dom}\bigg[x(k) - \gamma(k)\bigg(\nabla f\big(x(k)\big) + \frac{\partial g}{\partial x}\big(x(k)\big)^T\mu(k) + \alpha(k)x(k)\bigg)\bigg]
\end{equation}
and
\begin{equation} \label{eq:bigmu0}
\mu(k+1) = \bigg[\mu(k) + \gamma(k)\big(g\big(x(k)\big) - \alpha(k)\mu(k)\big)\bigg]_{+}.
\end{equation}

Above, $k$ is the iteration number, and $\alpha(k)$ and $\gamma(k)$ are defined as
\begin{equation}
\gamma(k) = \bar{\gamma}k^{-c_1} \textnormal{ and } \alpha(k) = \bar{\alpha}k^{-c_2}.
\end{equation}
The constants $\bar{\gamma}$, $\bar{\alpha}$, $c_1$, and $c_2$ are selected by the user,
with $c_1$ and $c_2$ subject to the conditions
\begin{equation}
0 < c_2 < c_1 \textnormal{ and } c_1 + c_2 < 1.
\end{equation}

One appealing aspect of the above method
in this setting is its robustness to noise in the values of $g$
and $\frac{\partial g}{\partial x}$. Define $G(x(k)) = g(x(k)) + w_g(k)$ and 
$\partial_x G(x(k)) = \frac{\partial g}{\partial x}\big(x(k)\big) + w_{\partial g}(k)$, where $w_g(k)$ and
$w_{\partial g}(k)$ are i.i.d. noises of the appropriate dimension. 
Then the noisy forms of Equations \eqref{eq:bigx0} and \eqref{eq:bigmu0} are
\begin{equation} \label{eq:bigx2}
x(k+1) = P_{\dom}\bigg[x(k) - \gamma(k)\big(\nabla f\big(x(k)\big) + \partial_{x}G(x(k))^T\mu(k) + \alpha(k)x(k)\big)\bigg]
\end{equation}
and
\begin{equation} \label{eq:bigmu2}
\mu(k+1) = \big[\mu(k) + \gamma(k)\big(G(x(k)) - \alpha(k)\mu(k)\big)\big]_{+}.
\end{equation}
It is this noisy form of update rule that will be used in the remainder of the paper. 
We use it to define Algorithm 1, below. 

\begin{samepage}
\emph{Algorithm 1} \\
\emph{Step 0:} Select a starting point $\big(x(0), \mu(0)\big) \in \dom \times \R^m$ and constants,
$\bar{\alpha}$, $\bar{\gamma}$, $c_1$, and $c_2$. Set $k = 0$. \\
\emph{Step 1:} Compute
\begin{equation} \label{eq:bigx}
x(k+1) = P_{\dom}\big[x(k) - \gamma(k)\big(\nabla f\big(x(k)\big)  + \partial_{x}G(x(k))^T\mu(k) + \alpha(k)x(k)\big)\big],
\end{equation}
\begin{equation} \label{eq:bigmu}
\mu(k+1) = \big[\mu(k) + \gamma(k)\big(G(x(k)) - \alpha(k)\mu(k)\big)\big]_{+}.
\end{equation}
\emph{Step 2:} Set $k = k+1$ and return to Step 1. \hfill $\triangle$ \\
\end{samepage}

For the purpose of analyzing the convergence
of Algorithm 1, we have the following definition.
\begin{definition} \label{def:conv}
Let $(\hat{x}, \hat{\mu})$ denote the unique saddle point of the Lagrangian as defined
in Equation \eqref{eq:bigl}. 
We say Algorithm 1 converges if it generates a sequence
$\big(x(k), \mu(k)\big)$ that converges in mean-square to $(\hat{x}, \hat{\mu})$
in the Euclidean norm, i.e.,
if
\begin{equation}
\lim_{k \to \infty} \mathbb{E}\big(\|x(k) - \hat{x}\|_2^2\big) = 0
\end{equation}
and
\begin{equation}
\lim_{k \to \infty} \mathbb{E}\big(\|\mu(k) - \hat{\mu}\|_2^2\big) = 0.
\end{equation}
\hfill $\triangle$
\end{definition} 

We now present the following theorem concerning the convergence
of Algorithm 1.

\begin{theorem} \label{thm:conv}
Algorithm 1 converges in the sense
of Definition \ref{def:conv} if the following four conditions are met:
\begin{enumerate}
\item $f_i$ and $g_j$ are convex and $C^1$ for all $i$ and $j$ \label{cond:convex}
\item $\mathcal{X}$ is a convex, closed, bounded set \label{cond:dom}
\item Slater's condition holds, i.e., there is some $\bar{x} \in \dom$ such that $g(\bar{x}) < 0$ \label{cond:slater}
\item all noises are zero mean, have bounded variance, and are independent at different points \label{cond:noise}
\end{enumerate}
\end{theorem}
\emph{Proof:} See \cite{poljak78}. \hfill $\blacksquare$

By assumption, Conditions \ref{cond:convex}, \ref{cond:dom}, and \ref{cond:slater} hold. It remains to
be shown that Condition \ref{cond:noise} can hold when differential privacy is desired. 
This will be shown in Section \ref{sec:main}. Below we 
describe an implementation of Algorithm 1 and then cover 
differential privacy. 


\subsection{Cloud Architecture}

Let the conditions and assumptions of Section \ref{ss:overview} hold. 
In \cite{hale14}, a cloud-based architecture was 
used and it was assumed that there was no
inter-agent communication. 
This assumption is kept in force
throughout this paper to enable 
privacy of agents' states.
In this architecture, each agent stores and manipulates
its own state $x_i(k) \in \R$. 
Similarly, the cloud
stores and updates $\mu^c(k)$, a vector 
of KT multipliers, as well as $x^c(k)$, the vector
of each agent's state; it does
not share $x^c(k)$ with the agents but instead 
only uses it to compute values of $\mu^c$. 

At each timestep $k$, the agents
receive information from the cloud, update their states, 
and then send their updated states to the cloud. At the same time
that the agents are updating their states, the cloud is computing
an update to $\mu^c$ which will be sent to the agents in the
next transmission the cloud sends. 
As mentioned above, the agents do not
talk to each other at all. The role of the cloud is to serve as a
trusted central aggregator for information so that computations
based upon sensitive, global information, namely $\frac{\partial g}{\partial x^c_i}$ and updates to $\mu^c$, 
can be computed and the results disseminated to the agents without
any agent ever directly knowing another agent's state. 

Before any optimization takes place, 
let agent $i$ be initialized with $f_i$,
$\bar{\alpha}$, $\bar{\gamma}$, $\dom$, $c_1$, $c_2$, and some initial state, $x_i(0)$. 
Let the cloud be initialized with $g$, $\frac{\partial g}{\partial x^c_i}$ for every $i$,
 $\bar{\alpha}$, $\bar{\gamma}$, $c_1$, $c_2$, and some
initial KT vector, $\mu^c(0)$, that is not based on the values of any initial states. 
To initialize the system, 
each agent sends its state to the cloud. Then
at each timestep, $k$, three actions occur.
First, agent $i$ receives a transmission containing 
private versions 
$\frac{\partial g}{\partial x^c_i}\big(x^c(k)\big)$
and $\mu^c(k)$ from the cloud (the details of the privacy
will be explained in Section \ref{sec:main}). Second, agent $i$ computes $x_i(k+1)$, and
simultaneously the cloud computes $\mu^c(k+1)$.
Third, agent $i$ sends $x_i(k+1)$ to the cloud, and
then this cycle of communication and computation is repeated with the newly computed values.


Writing out the update equations for a multi-agent implementation of 
Algorithm 1 using the cloud architecture, we see that
agent $i$ updates according to
\begin{equation} \label{eq:agentix}
x_i(k+1) = P_{\dom}\Bigg[x_i(k) - \gamma(k)\bigg(\frac{df_i}{dx_i}\big(x_i(k)\big) + \partial_{x_i}G(x^c(k))^T\mu^c(k) + \alpha(k)x_i(k)\bigg)\Bigg]
\end{equation}
and the cloud updates according to
\begin{equation} \label{eq:cloudmu}
\mu^c(k+1) = \bigg[\mu^c(k) + \gamma(k)\big(G(x^c(k)) - \alpha(k)\mu^c(k)\big)\bigg]_{+},
\end{equation}
where $\partial_{x_i} G(x^c(k)) = \frac{\partial g}{\partial x^c_i}\big(x^c(k)\big) + w_i(k)$ and,
as before, $G\big(x^c(k)\big) = g(x^c(k)) + w_g(k)$. 
The distributions of all noises $w_i$ and $w_g$ will be defined in Section \ref{sec:main}.
In Equation \eqref{eq:agentix}, the notation $\mu^c(k)$ is meant to indicate that at time $k$,
agent $i$ updates its state using the KT vector it just received from the cloud computer;
before the cloud computes a new value of $\mu^c$, this KT vector is equal to the $\mu$ vector stored on the cloud
at time $k$, $\mu^c(k)$, and hence is denoted as such. 

Due to the structure of communications in the system, agent $i$ receives $\partial_{x_i}G(x(k))$ and $\mu^c(k)$
before computing $x_i(k+1)$. Similarly, the cloud receives the agents' states at time $k$
(which end up being the contents of the vector $x^c(k)$) before computing $\mu^c(k+1)$. Then despite
the distributed nature of the problem, all information in the network is synchronized whenever updates
are computed anywhere. This means that, 
in aggregate, the steps taken by the agents and cloud using Equations \eqref{eq:agentix} and \eqref{eq:cloudmu} 
are identical to those used
in Algorithm 1 to solve
the global optimization problem defined in Section \ref{ss:overview}. 
As a result, the analysis
for the convergence of the cloud-based implementation of
Algorithm 1 will be carried out using 
the centralized form of Algorithm 1 as presented in Section \ref{ss:overview},
and it will apply to the cloud-based problem. 

\subsection{Differential Privacy for Dynamic Systems}

Let there be $n$ input signals to a system,
each contributed by some user. Let the $i^{th}$
input signal be denoted $u_i \in \ell^{s_i}_{p_i}$.
Here $s_i \in \N$ is the dimension of the signal
and, with an abuse of notation, we say $u_i \in \ell^{s_i}_{p_i}$ 
if each \emph{finite} truncation of $u_i$ has finite $p_i$-norm, i.e., if
\begin{equation}
u_{0:k} := (u(0)^T, u(1)^T, \ldots, u(k)^T)^T
\end{equation}
has finite $p_i$-norm for all $k$. 

Using this definition, the full input space to the system is
\begin{equation}
\ell_p^s = \ell^{s_1}_{p_1}\times\ell^{s_2}_{p_2}\times \cdots \times \ell^{s_n}_{p_n},
\end{equation}
and the system generates outputs $y \in \ell^r_{q}$ with $q, r \in \N$.
Fix a set of non-negative real numbers $b = (b_1, \ldots, b_n)$. We
define a symmetric binary adjacency relation, $\Adj(\cdot, \cdot)$, on the space $\ell^s_p$
such that 
\begin{equation}
\Adj_b(u, \tilde{u}) = 1 \textnormal{ if and only if } \|u_i - \tilde{u}_i\|_{p_i} \leq b_i\textnormal{, } \textnormal{for some $i$ and } u_j = \tilde{u}_j \textnormal{ for all } j \neq i.
\end{equation}
In words, $\Adj(u, \tilde{u})$ holds if and only if $u$ and $\tilde{u}$
differ by (at most) one component and this difference is bounded by the corresponding element of $b$. 
In this paper, we will focus
exclusively on the case that ${p_i = 2}$ for every $i$. 
The symbol $\|\cdot\|_2$ will be used
for both the Euclidean norm and the $\ell_2$ norm, though the meaning of each use
is clear from context.

Roughly, differential privacy guarantees that
if two input signals are adjacent, their output signals should not differ by too much.
As a result, small changes to inputs are not seen at the output and someone, e.g., a malicious eavesdropper,
with access to the output of the of the system cannot exactly determine the input. 
To make this notion precise, 
let us fix a probability space $(\Omega, \mathcal{F}, \P)$. 
Let $\mathcal{R}^d$ denote the Borel $\sigma$-algebra defined on $\mathbb{R}^d$. 
In the setting of differentially private dynamic systems, a \emph{mechanism} is a stochastic process $M$ of the form
\begin{equation}
M : \ell^{s}_{p} \times \Omega \to \ell^r_q,
\end{equation}
i.e., $M$ has inputs in $\ell^{s}_{p}$ and sample paths in $\ell^{r}_{q}$.
We now state a lemma concerning differentially private mechanisms for dynamical systems.
\begin{lemma} \label{lem:edprivacy}
A mechanism $M : \ell^s_p \times \Omega \to \ell^r_{q}$ is $(\epsilon,\delta)$-differentially private
if and only if for every $u, \tilde{u} \in \ell^s_p$ satisfying $\Adj_b(u, \tilde{u})$, we have
\begin{equation} \label{eq:ed}
\P\big((Mu)_{0:k} \in S\big) \leq e^{\epsilon}\P\big((M\tilde{u})_{0:k} \in S) + \delta, \forall k \geq 0, \forall S \in \mathcal{R}^{(k+1)r}.
\end{equation}
\end{lemma}
\emph{Proof:} See \cite{pappas14}, Lemma 2. \hfill $\blacksquare$

Equation \eqref{eq:ed} captures in a precise way the notion that, at each time $k$, 
truncated outputs up until $k$ that correspond to adjacent inputs must have similar
probability distributions, with the level of similarity of the outputs
determined by $\delta$ and $\epsilon$. 
Indeed, the constants $\epsilon$ and $\delta$ determine the level
of privacy afforded by the mechanism $M$ to users contributing input signals. 
Generally, $\epsilon$ is kept small, e.g., $0.1$, $\ln 2$, or $\ln 3$. The parameter $\delta$
should be kept very small because it allows a zero probability event for $\tilde{u}$ to
be a non-zero probability event for $u$ and thus controls
when there can noticeable differences in outputs which correspond to adjacent inputs.
The choice of $b$ determines which inputs to the system should produce similar outputs and
thus determines which inputs should ``look alike'' at the output. 

One appealing aspect of differential privacy is that its privacy guarantees cannot be weakened
by post-processing. Given the output of a differentially private mechanism, that output
can be processed freely without threatening the privacy of the inputs. 
In addition, this privacy guarantee holds
even against an adversary with arbitrary side information. Even if an adversary
gains knowledge of, e.g., some number of inputs, that adversary still cannot determine
exactly the system's other
input signals by observing its outputs. 
In the present paper, Equation \eqref{eq:ed} will be used as the definition
of differential privacy and the reader is referred to \cite{pappas14} 
for a proof of Lemma \ref{lem:edprivacy}.

Before discussing the mechanisms to be used here, we first define the $\ell_2$ sensitivity of
a system; while the $\ell_p$ sensitivity can be used for
other values of $p$, we focus on $p=2$ in this paper. Let $\mathcal{G}$ be a deterministic causal system. The 
$\ell_2$ sensitivity is an upper bound, $\Delta_2$, on the norm of the difference between the outputs
of $\mathcal{G}$
which correspond to adjacent inputs. That is, $\Delta_2$ must satisfy
\begin{equation}
\|\mathcal{G}u - \mathcal{G}\tilde{u}\|_2 \leq \Delta_2
\end{equation}
whenever $\Adj_b(u, \tilde{u})$ holds. 
There are several established differentially private mechanisms in the literature, e.g., \cite[Chapter 3]{dwork13},
though here we will use only the Gaussian mechanism in the lemma below. In it, we use the
$\Q$-function, defined as
\begin{equation}
\Q(x) = \frac{1}{\sqrt{2\pi}}\int_{x}^{\infty} e^{-\frac{u^2}{2}}du.
\end{equation}

\begin{lemma} \label{lem:gauss}
The mechanism $Mu = \mathcal{G}u + w$ is $(\epsilon, \delta)$-differentially private if
$w \sim \mathcal{N}(0, \sigma^2 I_r)$, where $I_r$ is the identity matrix of the same dimension as the output
space, $\ell_q^r$, and where $\sigma$ satisfies
\begin{equation}
\sigma \geq \Delta_2\kde
\end{equation}
where we define $K_{\delta} = \mathcal{Q}^{-1}(\delta)$ and
\begin{equation} \label{eq:kde}
\kde := \frac{1}{2\epsilon}\left(K_{\delta} + \sqrt{K_{\delta}^2 + 2\epsilon}\right).
\end{equation}
\end{lemma}
\emph{Proof:} See \cite{pappas14}. \hfill $\blacksquare$

We will use $\kde$ in Equation \eqref{eq:kde} for the remainder of the paper. 
Lemma $2$ says in a rigorous way that adding noise to the true output of a system can make
it differentially private, provided the distribution from which the noise is drawn
has large enough variance. 
Also from Lemma \ref{lem:gauss}, we see that once $\epsilon$ and $\delta$ are chosen, we need only
to find the $\ell_2$ sensitivity of a system, $\Delta_2$, to calibrate the level of noise
that must be added to guarantee differential privacy. We
do this in the next section for the optimization problem outlined above.

\section{Private Optimization} \label{sec:main}
Based on the optimization algorithm devised in Section \ref{sec:review}
and the cloud-based architecture, it is clear that two pieces of
information must be shared with agent $i$ at each time $k$: 
$\frac{\partial g}{\partial x^c_i}\big(x(k)\big)$ and $\mu^c(k)$. 
In order to keep the agents' states private from each other, we must
then alter these quantities, or the way in which
they are computed, before they are sent to the agents. 

We regard each $\frac{\partial g}{\partial x^c_i}$ as a dynamic system with output $y_i(k) \in \ell_2^m$
at time $k$, where $m$ is the number of constraints as defined in Section \ref{ss:overview}.
Noise is added to each $y_i$ before it is sent to the agents in order to guarantee differential privacy
of each agent's state. 
While agent $i$ will only make use of $\frac{\partial g}{\partial x^c_i}$
in its computations, we assume that each
agent has access to each output, i.e. agent $i$ has knowledge
of $\frac{\partial g}{\partial x^c_j}$ at each time even if $j \neq i$.

To ensure that no state values can be determined from $\mu^c$, we use the
resilience of differential privacy to post-processing. 
Specifically, rather than
adding noise to $\mu^c$ directly each time it is sent to the agents,
we regard $g$ as a dynamic system and add noise
to the value of $g$ when computing
$\mu^c$ in the cloud to guarantee that $g$ keeps $x^c$ differentially private. 
In doing this, $\mu^c$ then also keeps $x^c$ differentially private.
In regarding $g$ as a dynamical system and setting $y_g(k) = g(x^c(k))$, 
we see that its output $y_g \in \ell_2^m$, where $m$ is the
number of constraints as defined above.
In this framework, $\mu^c$ and each $\frac{\partial g}{\partial x^c_i}$ conceal each other
agent's state from agent $i$ and any eavesdropper intercepting the cloud's transmissions
to the agents. 

To define differentially private mechanisms for $\frac{\partial g}{\partial x^c_i}$, 
we first find bounds on the $\ell_2$ sensitivity of $\frac{\partial g}{\partial x^c_i}$. 
By assumption, $g \in C^2$ so that $\frac{\partial^2 g}{\partial x_i^{c,2}}$
is continuous. Also by assumption, $\dom$ is compact so that
$\left\|\frac{\partial^2 g}{\partial x_i^{c,2}}(x)\right\|_2$ attains
its maximum on $\dom$. Denote this maximum by $K^i_{g'}$. We see
that $K^i_{g'}$ is the Lipschitz constant for $\frac{\partial g}{\partial x^c_i}$.
Regarding $\frac{\partial g}{\partial x^c_i}$ as a system 
with inputs in $\R^n$ and outputs in $\R^m$, we have the following theorem
concerning a differentially private mechanism whose output will be
released to the agents. 
\begin{theorem} \label{thm:mech1}
Let $K_{g'}^i$ be the Lipschitz constant of $\frac{\partial g}{\partial x^c_i}$ over 
the domain $\dom$ and define
$B = \max_{i \in A}\{b_i\}$. Then the mechanism
\begin{equation} \label{eq:mech1mech}
M_{g'}^ix^c(k) = \frac{\partial g}{\partial x^c_i}\big(x(k)\big) + w_i(k)
\end{equation}
is differentially private
with $w_i(k) \sim \mathcal{N}(0, \sigma_i^2I_m)$, where $m$ is the number of constraints,
and with $\sigma_i$ satisfying
\begin{equation}
\sigma_i \geq \kde K_{g'}^iB.
\end{equation}
\end{theorem}
\emph{Proof:}
For two signals, $x$ and $\tilde{x}$ in $\ell_p^s$, satisfying
$\Adj_b(x, \tilde{x})$ for $b = (b_1, \ldots, b_n)$, we have
\begin{equation}
\left\|\frac{\partial g}{\partial x_i^c}(x(k)) - \frac{\partial g}{\partial x_i^c}(\tilde{x}(k))\right\|_2 \leq K^i_{g'}\|x(k) - \tilde{x}(k)\|_2 
\leq K^i_{g'}\|x - \tilde{x}\|_2 \leq K^i_{g'}b_i \leq K^i_{g'}B.
\end{equation}
Using this bound on the $\ell_2$ sensitivity of $\frac{\partial g}{\partial x^c_i}$, 
for each $i \in A$ we define $\sigma_i$ by
\begin{equation} \label{eq:sigmai}
\sigma_i \geq \kde K^i_{g'}B,
\end{equation}
with $\kde$ defined as before. The theorem is then a straightforward application of Lemma \ref{lem:gauss}.
\hfill $\blacksquare$



To keep the agents' states private in releasing $\mu$ we rely on the resilience to 
post-processing of differential privacy and compute $\mu$ using a differentially
private form of $g$. As above, we 
note that $\frac{\partial g}{\partial x}$ itself is continuous by assumption and
that there exists some $K_{g}$ satisfying
\begin{equation}
\left\|\frac{\partial g}{\partial x}(x)\right\|_2 \leq K_{g}
\end{equation}
for all $x \in \dom$. By definition, $K_{g}$ is the Lipschitz constant of $g$ and
we use this to define a mechanism for making $g$ private in the next theorem below.
We omit the proof of this theorem due to its similarity to 
the proof of
Theorem \ref{thm:mech1}
above. 

\begin{theorem} \label{thm:mech2}
Let $K_{g}$ be the Lipschitz constant of $g$ over $\dom$ and define
$B = \max_{i \in A}\{b_i\}$. Then the mechanism
\begin{equation} 
M_{g}x^c(k) = g\big(x^c(k)\big) + w_g(k)
\end{equation}
is differentially private
with $w_g(k) \sim \mathcal{N}(0, \sigma_g^2I_m)$, where $m$ is the number of constraints,
and with $\sigma_g$ satisfying
\begin{equation} \label{eq:sigmag}
\sigma_g \geq \kde K_{g}B.
\end{equation}
\end{theorem} \hfill $\blacksquare$


Under the assumptions already made in crafting the optimization problem in Section \ref{sec:review}, 
conditions $1-3$ of Theorem \ref{thm:conv} are satisfied. 
To satisfy condition $4$ of Theorem \ref{thm:conv}, we need only
to select values of $\sigma_i$ and $\sigma_g$ that are bounded above. Then Algorithm
1 will converge. 
In addition, as long as $\sigma_i$ and $\sigma_g$ are bounded below
as in
Equations \eqref{eq:sigmai} and \eqref{eq:sigmag}, respectively, the conditions
for $(\epsilon, \delta)$-differential privacy 
are simultaneously met.
Importantly,
Algorithm 1 is robust to noise appearing in exactly the places
it is injected for privacy. Adding more noise will increase the time
required for Algorithm 1
to converge and thus there is a natural trade-off: increased privacy
results in a slower convergence rate because more noise is added, while 
decreased privacy allows for faster convergence because it reduces the noise added.


\section{Simulation Results}
A numerical simulation was conducted to support the theoretical
developments of this paper. The problem here involves $n = 7$
agents and $m = 4$ constraints. The set $\dom = [-10, 10]^7$.
This set can represent, e.g., the area a team of robots
must stay inside in order to maintain communication links.
The constraint function $g$ was chosen to be
\begin{equation}
g(x) = \left(\begin{array}{c}
               x_1 + x_2 + x_3 - 3 \\
              x_5^2 + \frac{1}{12}x_6^4 + \frac{1}{12}x_7^4 - 20 \\
             x_3^2 + x_4 + x_6 - 1\\
            x_6^2 + x_7^2 - 5 \end{array}\right).
\end{equation}
Over $\dom$, the Lipschitz constant of $g$ was found to be approximately
$K_{g} = 472.567$. 
Each dynamic system $\frac{\partial g}{\partial x^c_i}$ and its Lipschitz constant over the domain
$\dom$ is listed in Table \ref{tab:dg}.

\begin{table}[!htbp]
\centering
\begin{tabular}{|M|M@{}@{\hspace{20pt}}|M|}
\hline
$i$ & $\frac{\partial g}{\partial x^c_i}$ & $K_{g'}^i$ \\[6pt] \hline \hline 
$1$ & $(1, 0, 0, 0)^T$ & $0$ \\[3pt] \hline
$2$ & $(1, 0, 0, 0)^T$ & $0$ \\[3pt] \hline
$3$ & $(1, 0, 2x_3, 0)^T$ & $2$ \\[3pt] \hline
$4$ & $(0, 0, 1, 0)^T$ & $0$ \\[3pt] \hline 
$5$ & $(0, 2x_5, 0, 0)^T$ & $2$ \\[3pt] \hline 
$6$ & $(0, \frac{1}{3}x_6^3, 1, 2x_6)^T$ & $100.08$ \\[3pt] \hline
$7$ & $(0, \frac{1}{3}x_7^3, 0, 2x_7)^T$ & $100.08$ \\[3pt] \hline
\end{tabular} 
\caption{Derivatives of $g$ and their Lipschitz constants}
\label{tab:dg}
\end{table}

Although the Lipschitz constants in Table \ref{tab:dg} are quite different,
it was desired to have the level of privacy guaranteed by each
mechanism be identical. The values $\delta = 0.05$ and
$\epsilon = \ln 3$ were chosen to be used
by all systems, giving ${\kde = 1.7565}$.
In addition, it was desired to make agent $i$'s state
difficult to distinguish from 
a ball of radius $1$ around it in the
space $\ell^{1}_2$, leading to the choice 
of $b_i = 1$ for all $i$, giving
\begin{equation}
b = (1, 1, 1, 1, 1, 1, 1)^T.
\end{equation}
Using these values, the distributions of noises to be added 
to each $\frac{\partial g}{\partial x^c_i}$ were
determined and are shown in Table \ref{tab:noise}.

\begin{table}[!htbp]
\centering
\begin{tabular}{|M|M@{}@{\hspace{30pt}}|}
\hline
Noise & Distribution \\[6pt] \hline \hline
$w_1$ & $\mathcal{N}(0, 0)$ \\[3pt] \hline
$w_2$ & $\mathcal{N}(0, 0)$ \\[3pt] \hline
$w_3$ & $\mathcal{N}(0, 12.3406I_{4})$ \\[3pt] \hline
$w_4$ & $\mathcal{N}(0, 0)$ \\[3pt] \hline
$w_5$ & $\mathcal{N}(0, 12.3406I_{4})$ \\[3pt] \hline 
$w_6$ & $\mathcal{N}(0, 30900.7580I_{4})$ \\[3pt] \hline 
$w_7$ & $\mathcal{N}(0, 30900.7580I_{4})$ \\[3pt] \hline
\end{tabular} 
\caption{Distributions of noise added to each $\frac{\partial g}{\partial x^c_i}$}
\label{tab:noise}
\end{table}
In addition, the noise added to $g$ when computing
$\mu$ was
\begin{equation}
w_g \sim \mathcal{N}(0, 688971.6017I_{4}),
\end{equation}
and the sum of the agents' objective functions was 
\begin{equation}
f(x) = (x_1 - 9)^2 + x_1 + (x_2 + 4)^4 + (x_3 - 1)^8 + x_4^2 + (x_4 + 6) + (x_5 + 3)^6 + (x_6 - 7)^2 + (x_7 - 5)^2. 
\end{equation}

For Algorithm 1, the values $\bar{\gamma} = 0.0005$,
$\bar{\alpha} = 0.20$, $c_1 = \frac{1}{3}$ and
$c_2 = \frac{1}{4}$
were chosen.
Having selected everything needed, a simulation was run that consisted
of $500,000$ timesteps. The primal and dual components 
of the saddle point of $L$ were computed ahead of
time to be approximately
\begin{equation}
\hat{x} = \left(7.591, -4.769, 0.178, -0.822, -2.863, 1.790, 1.340\right)^T
\end{equation}
and
\begin{equation}
\hat{\mu} = (1.8139, 0, 0.6409, 2.7314)^T. 
\end{equation}

\begin{figure}
\centering
\includegraphics[width=3.4in]{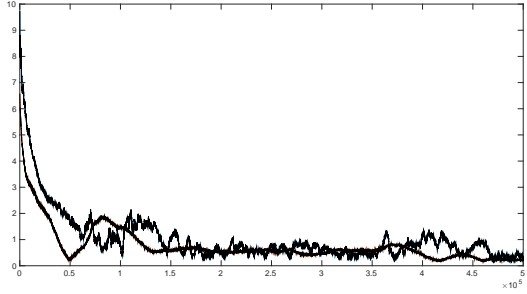}
\caption{Values of $\|x(k) - \hat{x}\|$ (jagged line) and $\|\mu(k) - \hat{\mu}\|$ (smooth line) for $k = 1, \ldots, 500,000$.}
\label{fig:norm}
\end{figure}

To illustrate the convergence of the algorithm, the values of $\|x(k) - \hat{x}\|_2$ and
$\|\mu(k) - \hat{\mu}\|_2$ for $k = 1, \ldots, 500,000$ are plotted in Figure \ref{fig:norm}.
We see that the
distance from $x(k)$ to $\hat{x}$ decreases almost monotonically 
in early iterations and oscillates more later on, though there is
a discernible decreasing trend while it oscillates. 
Similarly,
we see that 
the distance from $\mu(k)$ to $\hat{\mu}$ also follows a generally decreasing trend over time.
 The error values after $200,000$ iterations were
$\|x(200,000) - \hat{x}\|_2 = 0.4839$ and ${\|\mu(200,000) - \hat{\mu}\|_2 = 0.5459}$,
indicating close convergence to the unique saddle point after 
the first $200,000$ iterations.
Numerically, these distances to the optimum are generally within the acceptable
tolerance of error for many applications. For applications with tighter
bounds on the required distance to the optimum, the algorithm can simply be run for
longer. In this case, after $500,000$ iterations the error values were
$\|x(500,000) - \hat{x}\|_2 = 0.2612$ and ${\|\mu(500,000) - \hat{\mu}\|_2 = 0.2123}$,
indicating noticeable decreases in the level of sub-optimality in $x(k)$ and $\mu(k)$. 




The relative value of reducing oscillations and increasing privacy
will vary between problems and should be considered when
designing problems and selecting $\epsilon$, $\delta$, and $b$. Regardless of the weight of these two objectives,
differentially private multi-agent optimization successfully converges when agents
do not directly share their states, and indeed converges when it is impossible for any agent to discover
exactly any other agent's state value. 


\section{Conclusion}
A multi-agent optimization problem was considered in which it
is desirable to keep each agent's state private. This was achieved
via a primal-dual optimization method and
by using differential privacy. 
It was shown
that the conditions required for differential privacy and 
for convergence of the optimization algorithm can be simultaneously
satisfied and thus that differentially private multi-agent
optimization with constraints is possible. Numerical results were
then provided to attest to the viability of this approach.

\bibliographystyle{plain}{}
\bibliography{sources}

\end{document}